# Exponential Ratio–Type Estimators In Stratified Random Sampling

†Rajesh Singh, Mukesh Kumar, R. D. Singh, M. K. Chaudhary

Department of Statistics, B.H.U., Varanasi (U.P.)-India

†Corresponding author

---

**Abstract**

Kadilar and Cingi (2003) have introduced a family of estimators using auxiliary information in stratified random sampling. In this paper, we propose the ratio estimator for the estimation of population mean in the stratified random sampling by using the estimators in Bahl and Tuteja (1991) and Kadilar and Cingi (2003). Obtaining the mean square error (MSE) equations of the proposed estimators, we find theoretical conditions that the proposed estimators are more efficient than the other estimators. These theoretical findings are supported by a numerical example.

## 1. **Introduction**

Let a finite population having N distinct and identifiable units be divided into L strata. Let $n_h$ be the size of the sample drawn from $h^{th}$ stratum of size $N_h$ by using simple random sampling without replacement. Let

$$\sum_{h=1}^{L} n_h = n \quad \text{and} \quad \sum_{h=1}^{L} N_h = N.$$

Let y and x be the response and auxiliary variables respectively, assuming values $y_{hi}$ and $x_{hi}$ for the $i^{th}$ unit in the $h^{th}$ stratum.

Let the stratum means be

$$\overline{Y}_h = \frac{1}{N_h} \sum_{i=1}^{N_h} y_{hi} \text{ and } \overline{X}_h = \frac{1}{N_h} \sum_{i=1}^{N_h} x_{hi} \text{ respectively.}$$

A commonly used estimator for $\overline{Y}$ is the traditional combined ratio estimator defined as –

$$\overline{y}_{CR} = \frac{\overline{y}_{st}}{\overline{x}_{st}}.\overline{X}_{st} \tag{1.1}$$

where, $\overline{y}_{st} = \sum_{h=1}^{L} w_h \overline{y}_h$, $\quad \overline{x}_{st} = \sum_{h=1}^{L} w_h \overline{x}_h$,

$$\overline{y}_h = \frac{1}{n_h} \sum_{i=1}^{n_h} y_{hi} \quad \text{and} \quad \overline{x}_h = \frac{1}{n_h} \sum_{i=1}^{n_h} x_{hi},$$

$w_h = N_h/N$ and $\overline{X} = \sum_{h=1}^{L} w_h \overline{X}_h$.

The MSE of $\overline{y}_{CR}$, to a first degree of approximation, is given by

$$MSE(\overline{y}_{CR}) \cong \sum_{h=1}^{L} w_h^2 \gamma_h [S_{yh}^2 + R^2 S_{xh}^2 - 2RS_{yxh}] \tag{1.2}$$

where $\gamma_h = (\frac{1}{n_h} - \frac{1}{N_h})$, $R = \frac{\overline{Y}}{\overline{X}} = \frac{\overline{Y}_{st}}{\overline{X}_{st}}$ is the population ratio, $S_{yh}^2$ is the population variance of a variate of interest in stratum h and $S_{xh}^2$ is the population variance of auxiliary variate in stratum h and $S_{yxh}$ is the population covariance between auxiliary variate and variate of interest in stratum h.

Auxiliary variables are commonly used in survey sampling to improve the precision of estimates. Whenever there is auxiliary information available, the researchers want to utilize it in the method of estimation to obtain the most efficient estimator. In some cases, in addition to mean of auxiliary, various parameters related to auxiliary variable, such as standard deviation, coefficient of variation, skewness, kurtosis, etc. may also be known (Koyncu and Kadilar (2009). In recent years, a number of research papers on ratio type and regression type estimators have appeared, based on different types of transformation. Some of the contributions in this area are due to Sisodiya and Dwivedi (1981), Upadhyaya and Singh (1999), Singh and Tailor (2003), Kadilar and Cingi (2003, 2004, 2006), Singh et.al.(2004), Khoshnevisan et.al. (2007), Singh et.al. (2007) and Singh et.al. (2008). In this article, we study some of these transformations and propose an improved estimator.

## 2. Kadilar and Cingi Estimator

Kadilar and Cingi (2003) have suggested following modified estimator

$$\bar{y}_{KC} = \frac{\bar{y}_{st}}{\bar{x}_{st,ab}} \bar{X}_{st,ab} \tag{2.1}$$

where $\bar{x}_{st,ab} = \sum_{h=1}^{L} w_h \left( a\bar{x}_h + b \right)$, $\bar{X}_{st,ab} = \sum_{h=1}^{L} w_h \left( a\bar{X}_h + b \right)$.

and a, b suitably chosen scalars, these are either functions of the auxiliary variable x such as coefficient of variation $C_x$, coefficient of kurtosis $\beta_2(x)$ etc or some other constants.

The MSE of the estimator $\bar{y}_{KC}$ is given by

$$MSE(\bar{y}_{KC}) = \sum_{h=1}^{L} W_h^2 \gamma_h \left[S_{yh}^2 + R_{ab}^2 S_{xh}^2 - 2R_{ab} S_{yxh}\right] \tag{2.2}$$

where $R_{a,b} = \dfrac{\bar{Y}_{st}.\sum w_h \bar{X}_h.a_h}{\bar{X}_{st,ab}.(\bar{X}_{st})}$.

Bahl and Tuteja (1991) suggested an exponential ratio type estimator

$$\bar{y}_{BT} = \bar{y}\exp\left[\frac{\bar{X}-\bar{x}}{\bar{X}+\bar{x}}\right] \tag{2.3}$$

The estimator $\bar{y}_{BT}$ is more efficient than the usual ratio estimator under certain conditions. In recent years, many authors such as Singh et. al. (2007), Singh and Vishwakarma (2007) and Gupta and Shabbir (2007) have used Bahl and Tuteja (1991) estimator to propose improved estimators.

Following Bahl and Tuteja (1991) and Kadilar and Cingi (2003), we have proposed some exponential ratio type estimators in stratified random sampling.

## 3. **Proposed estimators**

The Bahl and Tuteja (1991) estimator in stratified sampling takes the following form

$$t = \bar{y}_{st}\exp\left[\frac{\bar{X}_{st}-\bar{x}_{st}}{\bar{X}_{st}+\bar{x}_{st}}\right] \tag{3.1}$$

The bias and MSE of t, to a first degree of approximation, are given by

$$Bias(t) = \frac{1}{\bar{X}_{st}}\sum_{h=1}^{L} w_h^2 \gamma_h \left(\frac{3R}{8}S_{xh}^2 - \frac{1}{2}S_{yxh}\right) \tag{3.2}$$

$$MSE(t) = \sum_{h=1}^{L} w_h^2 \gamma_h \left[S_{yh}^2 - RS_{yxh} + \frac{R^2}{4}S_{xh}^2\right] \tag{3.3}$$

3.1 **Sisodia- Dwivedi estimator**

When the population coefficient of variation $C_x$ is known, Sisodia and Dwivedi (1981) suggested a modified ratio estimator for $\bar{Y}$ as-

$$\bar{y}_{SD} = \bar{y}\frac{\bar{X}+C_x}{\bar{x}+C_x} \tag{3.4}$$

In stratified random sampling, using this transformation the estimator t will take the form

$$t_{SD} = \bar{y}\exp\left[\frac{\sum_{h=1}^{L} w_h(\bar{X}_h + C_{xh}) - \sum_{h=1}^{L} w_h(\bar{x}_h + C_{xh})}{\sum_{h=1}^{L} w_h(\bar{X}_h + C_{xh}) + \sum_{h=1}^{L} w_h(\bar{x}_h + C_{xh})}\right]$$

$$= \bar{y}_{st}\exp\left[\frac{\bar{X}_{SD} - \bar{x}_{SD}}{\bar{X}_{SD} + x_{SD}}\right] \tag{3.5}$$

where $\bar{X}_{SD} = \sum_{h=1}^{L} w_h(\bar{X}_h + C_{xh})$ and $\bar{x}_{SD} = \sum_{h=1}^{L} w_h(\bar{x}_h + C_{xh})$.

The bias and MSE of $t_{SD,}$ are respectively given by –

$$\text{Bias}(t_{SD}) = \frac{1}{\bar{X}_{st}}\sum_{h=1}^{L} w_h^2\gamma_h\theta_{SD}\left(\frac{R_{SD}}{8}S_{xh}^2 - \frac{1}{2}S_{yxh}\right) \tag{3.6}$$

$$\text{MSE}(t_{SD}) = \sum_{h=1}^{L} w_h^2\theta_h[S_{yh}^2 - R_{SD}S_{yxh} + \frac{R_{SD}^2}{4}S_{xh}^2] \tag{3.7}$$

where $R_{SD} = \frac{\bar{Y}_{st}}{\bar{X}_{SD}} = \frac{\sum_{h=1}^{L} w_h\bar{Y}_h}{\sum_{h=1}^{L} w_h(\bar{X}_h + C_{xh})}$ and $\theta_{SD} = \frac{\bar{X}_{st}}{\bar{X}_{SD}}$ .

3.2 **Singh-Kakran estimator**

Motivated by Sisodiya and Dwivedi (1981), Singh and Kakran (1993) suggested another ratio-type estimator for estimating $\overline{Y}$ as-

$$\overline{y}_{SK} = \overline{y}\frac{\overline{X}+\beta_2(x)}{\overline{x}+\beta_2(x)} \tag{3.8}$$

Using (3.8), the estimator t at (3.1) will take the following form in stratified random sampling-

$$t_{SK} = \overline{y}\exp\left[\frac{\sum_{h=1}^{L} w_h(\overline{X}_h + \beta_{2h}(x)) - \sum_{h=1}^{L} w_h(\overline{x}_h + \beta_{2h}(x))}{\sum_{h=1}^{L} w_h(\overline{X}_h + \beta_{2h}(x)) + \sum_{h=1}^{L} w_h(\overline{x}_h + \beta_{2h}(x))}\right]$$

$$= \overline{y}_{st}\exp\left[\frac{\overline{X}_{SK} - \overline{x}_{SK}}{\overline{X}_{SK} + \overline{x}_{SK}}\right] \tag{3.9}$$

where,

$$\overline{x}_{SK} = \sum_{h=1}^{L} w_h(\overline{x}_h + \beta_{2h}(x)) \text{ and } \overline{X}_{SK} = \sum_{h=1}^{L} w_h(\overline{X}_h + \beta_{2h}(x)).$$

Bias and MSE of $t_{SK}$, are respectively given by

$$\text{Bias}(t_{SK}) = \frac{1}{\overline{X}_{st}}\sum_{h=1}^{L} w_h^2\gamma_h\ \theta_{SK}\left(\frac{R_{SK}}{8}S_{xh}^2 - \frac{1}{2}S_{yxh}\right) \tag{3.10}$$

$$\text{MSE}(t_{SK}) = \sum_{h=1}^{L} w_h^2\gamma_h\left[S_{yh}^2 - R_{SK}S_{yxh} + \frac{R_{SK}^2}{4}S_{xh}^2\right] \tag{3.11}$$

where $R_{SK} = \dfrac{\overline{Y}_{st}}{\sum_{h=1}^{L} w_h(\overline{X}_h + \beta_{2h}(x))}$ and $\theta_{SK} = \dfrac{\overline{X}_{st}}{\overline{X}_{SK}}$.

3.3 **Upadhyaya-Singh estimator**

Upadhyaya and Singh (1999) considered both coefficients of variation and kurtosis in their ratio type estimator as

$$\bar{y}_{US1} = \bar{y}\,\frac{\bar{X}\beta_2(x) + C_x}{\bar{x}\beta_2(x) + C_x} \tag{3.12}$$

We adopt this modification in the estimator t proposed at (3.1)

$$t_{US1} = \bar{y}_{st}\exp\left(\frac{\sum_{h=1}^{L} w_h(\bar{X}_h\beta_{2h}(x) + C_{xh})}{\sum_{h=1}^{L} w_h(\bar{x}_h\beta_{2h}(x) + C_{xh}}\right) \tag{3.13}$$

$t_{US1}$ at (3.13) can be re-written as

$$t_{US1} = \bar{y}_{st}\exp\left[\frac{\bar{X}_{US1} - \bar{x}_{US1}}{\bar{X}_{US1} + \bar{x}_{US1}}\right] \tag{3.14}$$

where $\bar{X}_{US1} = \sum_{h=1}^{L} w_h(\bar{X}_h\beta_{2h}(x) + C_{xh})$ $\quad\quad$ $\bar{x}_{US1} = \sum_{h=1}^{L} w_h(\bar{x}_h\beta_{2h} + C_{xh})$.

Bias and MSE of $t_{US1}$, to first degree of approximation, are respectively given by

$$Bias(t_{US1}) = \frac{1}{\bar{X}_{st}}\sum_{h=1}^{L} w_h^2\gamma_h\theta_{US1}\left(\frac{R_{US1}}{8}S_{xh}^2 - \frac{1}{2}S_{yxh}\right) \tag{3.15}$$

$$MSE(t_{US1}) = \sum_{h=1}^{L} w_h^2\gamma_h\left[S_{yh}^2 - R_{US1}S_{yxh} + \frac{R_{US1}^2}{4}S_{xh}^2\right] \tag{3.16}$$

where $R_{us1} = \dfrac{\bar{Y}_{st}.\sum w_h\bar{X}_h.\beta_{2h}(x)}{\sum_{h=1}^{L} w_h(\bar{X}_h\beta_{2h}(x) + C_{xh})\bar{X}_{st}}$ and $\theta_{US1} = \dfrac{\sum w_h\bar{X}_h.\beta_{2h}(x)}{\bar{X}_{US1}}$.

Upadhyaya and Singh (1999) proposed another estimator by changing the place of coefficient of kurtosis and coefficient of variation as

$$\bar{y}_{US2} = \bar{y}\frac{\bar{X}C_x + \beta_2(x)}{\bar{x}C_x + \beta_2(x)} \qquad (3.17)$$

Incorporating this modification in the proposed estimator t, we have-

$$t_{US2} = \bar{y}_{st}\exp\left[\frac{\bar{X}_{US2} - \bar{x}_{US2}}{\bar{X}_{US2} + \bar{x}_{US2}}\right] \qquad (3.18)$$

where $\bar{x}_{us2} = \sum_{h=1}^{L} w_h(\bar{X}_h C_{xh} + \beta_{2h}(x))$ and $\bar{X}_{us2} = \sum_{h=1}^{L} w_h(\bar{X}_h C_{xh} + \beta_{2h}(x))$

Bias and MSE of $t_{US2,}$ are respectively given by –

$$Bias(t_{US2}) = \frac{1}{\bar{X}_{st}}\sum_{h=1}^{L} w_h^2 \gamma_h \theta_{US2}\left(\frac{R_{US2}}{8}S_{xh}^2 - \frac{1}{2}S_{yxh}\right) \qquad (3.19)$$

$$MSE(t_{US2}) = \sum_{h=1}^{L} w_h^2 \gamma_h [S_{yh}^2 - R_{US2}S_{yxh} + \frac{R_{US2}^2}{4}S_{xh}^2] \qquad (3.20)$$

where $R_{US2} = \dfrac{\bar{Y}_{st}.\sum w_h \bar{X}_h.C_{xh}}{\sum_{h=1}^{L} w_h\left(\bar{X}_h C_{xh} + \beta_{2h}(x)\right)\bar{X}_{st}}$

and $\theta_{US2} = \dfrac{\sum w_h \bar{X}_h C_{xh}}{\bar{X}_{US2}}$.

3.4. **G.N. Singh Estimator**

Following Singh (2001), using values of $\sigma_x$ and $\beta_2(x)$, we propose following two estimators.

$$t_{GNS1} = \bar{y}_{st}\exp\left[\frac{\bar{X}_{GNS1} - \bar{x}_{GNS1}}{\bar{X}_{GNS1} + \bar{x}_{GNS1}}\right] \qquad (3.21)$$

where $\overline{X}_{GNS1} = \sum_{h=1}^{L} w_h (\overline{X}_h + \sigma_{xh})$ and $\overline{x}_{GNS1} = \sum_{h=1}^{L} w_h (\overline{x}_h + \sigma_{xh})$

The Bias and MSE of $t_{GNS1}$ to a first degree of approximation, are respectively given by –

$$Bias(t_{GNS1}) \frac{1}{\overline{X}_{st}} \sum_{h=1}^{L} w_h^2 \gamma_h \theta_{GNS1} \left(\frac{R_{GNS1}}{8} S_{xh}^2 - \frac{1}{2} S_{yxh}\right) \tag{3.22}$$

$$MSE(t_{GNS1}) = \sum_{h=1}^{L} w_h^2 \gamma_h [S_{yh}^2 - R_{GNS1} S_{yxh} + \frac{R_{GNS1}^2}{4} S_{xh}^2] \tag{3.23}$$

where $R_{GNS1} = \dfrac{\overline{Y}_{st}}{\sum_{h=1}^{L} w_h \left(\overline{X}_h + \sigma_{xh}\right)}$, $\theta_{GNS1} = \dfrac{X_{st}}{\overline{X}_{GNS1}}$.

Similarly, we propose another estimator

$$t_{GNS2} = \overline{y}_{st} \exp\left[\frac{\overline{X}_{GNS2} - \overline{x}_{GNS2}}{\overline{X}_{GNS2} + x_{GNS2}}\right] \tag{3.24}$$

where, $\overline{X}_{GNS2} = \sum_{h=1}^{L} w_h (\overline{X}_h \beta_{2h}(x) + \sigma_{xh})$, $\overline{x}_{GNS2} = \sum_{h=1}^{L} w_h (\overline{x}_h \beta_{2h} + \sigma_{xh})$

The Bias and MSE of $t_{GNS2}$ to a first degree of approximation are respectively given by –

$$Bias(t_{GNS2}) = \frac{1}{\overline{X}_{st}} \sum_{h=1}^{L} w_h^2 \gamma_h \theta_{GNS2} \left(\frac{3R_{GNS2}}{8} S_{xh}^2 - \frac{1}{2} S_{yxh}\right) \tag{3.25}$$

$$MSE(t_{GNS2}) = \sum_{h=1}^{L} w_h^2 \gamma_h [S_{yh}^2 - R_{GNS2} S_{yxh} + \frac{R_{GNS2}^2}{4} S_{xh}^2] \tag{3.26}$$

where, $R_{GNS1} = \dfrac{\overline{Y}_{st}}{\sum_{h=1}^{L} w_h \left(\overline{X}_h \beta_{2h}(x) + \sigma_{xh}\right)}$, $\theta_{GNS2} = \dfrac{\sum_{i=1}^{L} w_h \overline{X}_h . \beta_{2h}(x)}{\overline{X}_{GNS2}}$

4. **Improved Estimator**

Motivated by Singh et. al. (2008), we propose a new family of estimators given by-

$$t_{MK} = \bar{y}_{st} \exp\left[\frac{\bar{X}_{st,ab} - \bar{x}_{st,ab}}{\bar{X}_{st,ab} + \bar{x}_{st,ab}}\right]^{\alpha} \tag{4.1}$$

where a and b are suitably chosen scalars and $\alpha$ is a constant.

The bias and MSE of $t_{MK}$ up to first order of approximation, are respectively given by

$$\text{Bias}(t_{MK}) = \frac{1}{\bar{X}_{st}} \sum_{h=1}^{L} w_h^2 \gamma_h \theta_{ab} \left( \frac{\alpha(\alpha+2) R_{ab}}{8} S_{xh}^2 - \frac{1}{2} S_{yxh} \right) \tag{4.2}$$

$$\text{MSE}(t_{MK}) = \sum_{h=1}^{L} w_h^2 \gamma_h \left[ S_{yh}^2 - \alpha R_{ab} S_{yxh} + \frac{\alpha^2 R_{ab}^2}{4} S_{xh}^2 \right] \tag{4.3}$$

The $\text{MSE}(t_{MK})$ is minimized for the optimal value of $\alpha$ given by-

$$\alpha = 2 \frac{\sum_{i=1}^{L} w_h^2 \gamma_h S_{yxh}}{R_{ab} \sum_{i=1}^{L} w_h^2 \gamma_h S_{xh}^2}$$

Putting this value of $\alpha$ in equation (4.3), we get the minimum MSE of the estimator $t_{MK}$ as-

$$\text{MSE}(t_{MK})_{min.} = \sum_{i=1}^{L} w_h^2 \gamma_h S_{yh}^2 (1 - \rho_c^2) \tag{4.4}$$

where, $\rho_c$ is combined correlation coefficient in stratified sampling across all strata. It is calculated as $\rho_c^2 = \frac{\left( \sum_{i=1}^{L} w_h^2 \gamma_h \rho_h S_{yh} S_{xh} \right)^2}{\sum_{1=1}^{L} w_h^2 \gamma_h S_{yh}^2 \sum_{i=1}^{L} w_h^2 \gamma_h S_{xh}^2}$.

We note here that min MSE of $t_{MK}$ is independent of a and b. therefore, we conclude that it same for any (all) values of a and b.

5. **Efficiency comparisons**

First we compare the efficiency of the estimator t at (3.1) with estimator $t_{SD.}$ We have

$MSE(t_{SD}) < MSE(t)$

$$\sum_{h=1}^{L} w_h^2 \gamma_h \left[ S_{yh}^2 - R_{SD} S_{yxh} + \frac{R_{SD}^2}{4} S_{xh}^2 \right] < \sum_{h=1}^{L} w_h^2 \gamma_h \left[ S_{yh}^2 RS_{yxh} + \frac{R^2}{4} S_{xh}^2 \right] \quad (5.1)$$

$$\sum_{h=1}^{L} w_h^2 \gamma_h \left[ -R_{SD} S_{yxh} + \frac{R_{SD}^2}{4} S_{xh}^2 \right] < \sum_{h=1}^{L} w_h^2 \gamma_h \left[ -RS_{yxh} + \frac{R^2}{4} S_{xh}^2 \right]$$

Let $A = \sum_{h=1}^{L} w_h^2 \gamma_h S_{yxh}$ and $B = \sum_{h=1}^{L} w_h^2 \gamma_h S_{xh}^2$

Then equation (5.1) can be re-written as-

$$-R_{SD} A + \frac{R_{SD}^2}{4}.B \quad < -RA + \frac{R^2}{4}.B$$

$$-A\,(R_{SD} - R) + B/4(R_{SD} - R)\,(R_{SD} + R) < 0 \quad (5.2)$$

From (5.2), we get two conditions

$$\sum_{h=1}^{L} w_h^2 \theta_h \left[ -R_{SD} S_{yxh} + \frac{R_{SD}^2}{4} S_{xh}^2 \right] < \sum_{h=1}^{L} w_h^2 \theta_h \left[ -RS_{yxh} + \frac{R^2}{4} S_{xh}^2 \right]$$

(i) When $(R_{SD}-R)\,(R_{SD}+R) > 0$

$$B < 4A\,/(R_{SD}+R) \quad (5.3)$$

(ii) When $(R_{SD} - R)\,(R_{SD} + R) < 0$

$$B > 4A\,/\,(R_{SD}+R) \quad (5.4)$$

where $A = \sum_{h=1}^{L} w_h^2 \theta_h S_{yxh}$ and $B = \sum_{h=1}^{L} w_h^2 \theta_h S_{xh}^2$

When either of these conditions is satisfied, estimator $t_{SD}$ will be more efficient than the estimator t.

The same conditions also holds true for the estimators $t_{SK}$, $t_{US1}$, $t_{US2}$, $t_{GNS1}$ and $t_{GNS2}$ if we replace $R_{SD}$ by $R_{SK}$, $R_{US1}$, $R_{US2}$, $R_{GNS1}$ and $R_{GNS2}$ respectively in conditions (i) and (ii).

Next we compare the efficiencies of $t_{opt}$ with the other proposed estimators.

$$MSE(t_{MK})_{min.} < MSE(t_{ab})$$

$$\sum_{i=1}^{L} w_h^2 \gamma_h S_{yh}^2 (1-\rho_c^2) < \sum_{i=1}^{L} w_h^2 \gamma_h \left[ S_h^2 + \frac{R_{ab}^2}{4} S_x^2 - R_{ab} S_{yx} \right] \qquad (5.5)$$

On putting the value of $\rho_c$ and rearranging the terms we get

$$\left( \sum_{i=1}^{L} w_h^2 \gamma_h S_{yxh} - \frac{R_{ab}}{2} \sum_{i=1}^{L} w_h^2 \gamma_h S_{xh}^2 \right)^2 > 0 \qquad (5.6)$$

This is always true. Hence the estimator $t_{MK}$ under optimum condition will be more efficient than other proposed estimators in all conditions.

## 6. **Data description and results**

For empirical study we use the data set earlier used by Kadilar and Cingi (2003). Y is apple production amount in 854 villages of turkey in 1999, and x is the numbers of apple trees in 854 villages of turkey in 1999.

The data are stratified by the region of turkey from each stratum, and villages are selected randomly using the Neyman allocation as

$$n_h = \frac{N_h S_h}{\sum_{h=1}^{L} N_h S_h}$$

**Table 6.1: Data Statistics**

| | | |
|---|---|---|
| $N_1$=106 | $N_2$=106 | $N_3$=94 |
| $N_4$=171 | $N_5$=204 | $N_6$=173 |
| $n_1$=9 | $n_2$=17 | $n_3$=38 |
| $n_4$=67 | $n_5$=7 | $n_6$=2 |
| $\overline{X}_1 = 24375$ | $\overline{X}_2 = 27421$ | $\overline{X}_3 = 72409$ |
| $\overline{X}_4 = 74365$ | $\overline{X}_5 = 26441$ | $\overline{X}_6 = 9844$ |
| $\overline{Y}_1 = i536$ | $\overline{Y}_2 = 2212$ | $\overline{Y}_3 = 9384$ |
| $\overline{Y}_4 = 5588$ | $\overline{Y}_5 = 967$ | $\overline{Y}_6 = 404$ |
| $\beta_{x1} = 25.71$ | $\beta_{x2} = 34.57$ | $\beta_{x3} = 26.14$ |
| $\beta_{x4} = 97.60$ | $\beta_{x5} = 27.47$ | $\beta_{x6} = 28.10$ |
| $C_{x1}$=2.02 | $C_{x2}$=2.10 | $C_{x3}$=2.22 |
| $C_{x4}$=3.84 | $C_{x5}$=1.72 | $C_{x6}$=1.91 |
| $C_{y1}$=4.18 | $C_{y2}$=5.22 | $C_{y3}$=3.19 |
| $C_{y4}$=5.13 | $C_{y5}$=2.47 | $C_{y6}$=2.34 |
| $S_{x1}$=49189 | $S_{x2}$=57461 | $S_{x3}$=160757 |
| $S_{x4}$=285603 | $S_{x5}$=45403 | $S_{x6}$=18794 |
| $S_{y1}$=6425 | $S_{y2}$=11552 | $S_{y3}$=29907 |
| $S_{y4}$=28643 | $S_{y5}$=2390 | $S_{y6}$=946 |
| $\rho_1 = 0.82$ | $\rho_2 = 0.86$ | $\rho_3 = 0.90$ |
| $\rho_4 = 0.99$ | $\rho_5 = 0.71$ | $\rho_6 = 0.89$ |
| $\gamma_1 = 0.102$ | $\gamma_2 = 0.049$ | $\gamma_3 = 0.016$ |
| $\gamma_4 = 0.009$ | $\gamma_5 = 0.138$ | $\gamma_6 = 0.006$ |
| $w_1^2 = 0.015$ | $w_2^2 = 0.015$ | $w_3^2 = 0.012$ |
| $w_4^2 = 0.04$ | $w_5^2 = 0.057$ | $w_6^2 = 0.041$ |
| N=854 | n=140 | $\beta_x = 312.07$ |
| $C_x$=3.85 | $C_y$=5.84 | $S_x$=144794 |
| $S_y$=17106 | $\rho = 0.92$ | $\overline{X} = 37600$ |
| $\overline{Y} = 2930$ | R =0.07793 | $R_{SD}$=0.07792 |
| $R_{SK}$=0.07784 | $R_{US1}$=.07789 | $R_{US2}$=0.07786 |
| $R_{GNS1}$= 0.06632 | $R_{GNS2}$=0.07765 | |

**Table 6.2 : Estimators with their MSE values**

| **Estimators** | **MSE values** |
|---|---|
| $t$ | 359619.594 |
| $t_{SD}$ | 359649.688 |
| $t_{SK}$ | 359890.313 |
| $t_{US1}$ | 359739.875 |
| $t_{US2}$ | 359830.125 |
| $t_{GNS1}$ | 360007.985 |
| $t_{GNS2}$ | 360479.192 |
| $t_{MK(opt)}$ | 218374.8898 |

From Table 6.2, we conclude that the estimator $t_{MK}$ has the minimum MSE and hence it is most efficient among the discussed estimators.

## 7. **Conclusion**

In the present paper we have examined the properties of exponential ratio type estimators in stratified random sampling. We have derived the MSE of the proposed estimators and also that of some modified estimators and compared their efficiencies theoretically and empirically.